\documentclass[11pt]{amsart}
\usepackage[all]{xy}
\usepackage{amsmath}
\usepackage{amsfonts}
\usepackage{amssymb}
\usepackage{latexsym}
\usepackage{graphicx}
\usepackage{color}
\usepackage{amscd}
\usepackage{epsfig}
\usepackage{cite}
\newtheorem{problem}{Problem}[section]
\newtheorem{definition}[problem]{Definition}
\newtheorem{lemma}[problem]{Lemma}
\newtheorem{theorem}[problem]{Theorem}

\newtheorem{corollary}[problem]{Corollary}

\title{Generic Twisted $T$-adic exponential sums of binomials}
\author{Chunlei Liu}
\address{Department of Mathematics, Shanghai Jiao Tong
University, Shanghai 200240, P.R. China, E-mail: clliu@sjtu.edu.cn}
\author{Chuanze Niu}\address{School of
Mathematical Sciences, Beijing Normal University, Beijing 100875,
P.R. China, E-mail: niuchuanze@mail.bnu.edu.cn}
\begin{document}
\maketitle
%==================================================================
\begin{abstract}
The twisted $T$-adic exponential sum associated to $x^{d}+\lambda x$ is studied.
If $\lambda\neq0,$ then
an explicit arithmetic polygon is proved to
be the Newton polygon of the twisted $C$-function of the
T-adic exponential sum. It gives the Newton polygons of the
$L$-functions of twisted $p$-power order exponential sums.
\end{abstract}

%{\it Key words}: exponential sum,
%$L$-function, deformation, Newton polygon

%{\it MSC2000}: 11L07, 14F30
%=================================================================

%=========================================================================================================
\section{Introduction to twisted exponential sums}

Let $W$ be a Witt ring scheme of Witt vectors, $\mathbb{F}_q$ the field
of characteristic $p$ with $q$ elements,
$\mathbb{Z}_{q}=W(\mathbb{F}_q)$, and
$\mathbb{Q}_q=\mathbb{Z}_q[\frac{1}{p}]$.

Let $\triangle\supsetneq\{0\}$ be an integral convex polytope
 in $\mathbb{R}^n$, and $I$ the set of vertices of
 $\triangle$ different from the origin.
Let
$$f(x)=\sum\limits_{v\in \triangle}(a_vx^v,0,0,\cdots)\in
W(\mathbb{F}_q[x_1^{\pm1},\cdots,x_n^{\pm1}])\text{ with }
\prod_{v\in I}a_v\neq0,$$ where $x^v=x_1^{v_1}\cdots x_n^{v_n}$ if
$v=(v_1,\cdots,v_n)\in\mathbb{Z}^n$.

Let $T$ be a variable. Let $\mu_{q-1}$ be the group of $q-1$-th roots of unity in $\mathbb{Z}_q$ and
 $\chi=\omega^{-u}$ with $u\in\mathbb{Z}^n/(q-1)$ a fixed multiplicative
character of $(\mathbb{F}_q^{\times})^n$ into $\mu_{q-1},$  where
$\omega: x\rightarrow \hat{x}$ is the Teichm\"{u}ller character.

\begin{definition}
The sum
$$S_{f,u}(k,T)=\sum_{x\in\mathbb{F}_{q^k}^{\times}}\chi({\rm Norm}_{\mathbb{F}_{q^k}/\mathbb{F}_q}(x))(1+T)^{Tr_{\mathbb{Q}_{q^k}/\mathbb{Q}_p}(\hat{f}(\hat{x}))}\in\mathbb{Z}_q[[T]]$$
is called a twisted $T$-adic exponential sum. And the function
$$L_{f,u}(s,T)=\exp(\sum_{k=1}^{\infty}S_{f,u}(k,T)\frac{s^k}{k})\in 1+s\mathbb{Z}_q[[T]][[s]]$$
is called an $L$-function of a twisted T-adic exponential sum.
\end{definition}

\begin{definition}
The function
$$C_{f,u}(s,T)=\exp(\sum_{k=1}^{\infty}-(q^k-1)^{-1}S_{f,u}(k,T)\frac{s^k}{k}),$$ is called a
 $C$-function of a twisted $T$-adic exponential sums.
\end{definition}

The L-function and C-function determine each other :
$$L_{f,u}(s,T)=\prod_{i=0}^n C_{f,u}(q^is,T)^{(-1)^{n-i+1}{n\choose i}},$$
and $$C_{f,u}(s,T)= \prod_{j=0}^{\infty} L_{f,u}(q^js,T)^{(-1)^{n-1}{n-j+1\choose j}}.$$
By the last identity, one sees that
$$C_{f,u}(s,T)\in 1+s\mathbb{Z}_q[[T]][[s]].$$

The $T$-adic exponential sums were first introduced by Liu-Wan \cite{LWn}.
We view $L_{f,u}(s, T)$ and $C_{f,u}(s, T)$ as power
series in the single variable $s$ with coefficients in the $T$-adic
complete field $\mathbb{Q}_q((T))$. The $C$-function $C_{f,u}(s,T)$ was
shown to be $T$-adic entire in $s$ by Liu-Wan \cite{LWn} for $u=0$ and
Liu \cite{Liu1} for all $u$.

Let $\zeta_{p^m}$ be a primitive $p^m$-th root of unity, and $\pi_m=\zeta_{p^m}-1.$ Then
$L_{f,u}(s,\pi_m)$ is the $L$-function of the $p$-power order
exponential sums $S_{f,u}(k,\pi_m)$ studied by Adolphson-Sperber \cite{AS1,AS2,AS3,AS4} for $m=1,$
by Liu-Wei \cite{LW} and Liu \cite{Liu1} for $m\geq1.$

Let $C(\triangle)$ be the cone generated by $\triangle$. There is a degree
function $\deg$ on $C(\triangle)$ which is
$\mathbb{R}_{\geq0}$-linear and takes values $1$ on each
co-dimension $1$ face not containing $0$. For $a\not\in
C(\triangle)$, we define $\deg(a)=+\infty$. Write $C_u(\triangle)=C(\triangle)\cap (u+(q-1)\mathbb{Z}^n)$ and
$M_u(\Delta)=\frac{1}{q-1}C_u(\triangle).$ Let $b$ be the least positive integer such that $p^bu\equiv u(\mod q-1).$
Order elements of $\bigcup_{i=0}^{b-1}M_{p^iu}(\Delta)$ so that $\deg(x_1)\leq\deg(x_2)\leq\cdots.$

\begin{definition}
The infinite u-twisted Hodge polygon $H_{\Delta,u}^{\infty}$ of $\Delta$ is the convex
function on $[0,+\infty]$ with initial point 0 which is linear between consecutive integers and whose slopes are
$$\frac{\deg(x_{bi+1})+\deg(x_{bi+2})+\cdots+\deg(x_{b(i+1)})}{b}, i=0,1,\cdots.$$
\end{definition}

Write NP for the short of Newton polygon.
Liu \cite{Liu2} proved the following Hodge bound for the $C$-function $C_{f,u}(s,T).$
\begin{theorem} We have
$$T-\text{adic NP of }C_{f,u}(s,T)\geq\text{ord}_p(q)(p-1)H_{\Delta,u}^{\infty}.$$
\end{theorem}

In the rest of this paper, $f(x)=(x^d,0,0,\cdots)+(\lambda x,0,0,\cdots)\in W(\mathbb{F}_q[x])$ where
 $\lambda\in\mathbb{F}_q^{\times}$ and $\Delta=[0,d].$ We fix $0\leq u\leq q-1.$
We shall study the twisted exponential sum of $f(x).$

Let $a=\log_pq.$ Write $u=u_0+u_1p+\cdots+u_{a-1}p^{a-1}$ with $0\leq u_i\leq p-1.$
Then we have
$$\frac{u}{q-1}=-(u_0+u_1p+\cdots),~~~u_i=u_{b+i}.$$
Write $p^iu=q_i(q-1)+s_i$ for $i\in \mathbb{N}$ with $0\leq
s_i<q-1,$ then $s_{a-l}=u_l+u_{l+1}p+\cdots+u_{a+l-1}p^{a-1}$ for
$0\leq l\leq a-1$ and $s_i=s_{b+i}.$

\begin{lemma}
The infinite u-twisted Hodge polygon $H_{[0,d],u}^{\infty}$ is the convex
function on $[0,+\infty]$ with initial point 0 which is linear between consecutive integers and whose slopes are
$$\frac{u_0+u_1+\cdots+u_{b-1}}{bd(p-1)}+\frac{l}{d},~~l=0,1,\cdots.$$
\end{lemma}
\proof Observe that $\frac{s_i}{q-1}+l$ with $0\leq i\leq b-1$ is just
a permutation of $x_{bl+1},x_{bl+2},\cdots,x_{b(l+1)}.$ The lemma
follows.
\endproof

\begin{definition}
The arithmetic polygon $P_{\{1,d\},u}$ is the convex function on $\mathbb{N}$ which
is linear between consecutive integers
with initial point 0 and whose slopes are
$$\omega(n)=\frac{(p-1)n+\frac{1}{b}\sum\limits_{i=0}^{b-1}u_i}{d}+(d-1)\frac{1}{b}\sum_{i=0}^{b-1}(\{\frac{pn+u_i}{d}\}-\{\frac{n}{d}\}),$$
where $\{\cdot\}$ is the fractional part of a real number.
\end{definition}
The  main results of this paper are the following.
\begin{theorem}\label{main of example} We have
$$P_{\{1,d\},u}\geq (p-1)H_{[0,d],u}^{\infty}.$$
Moreover, they coincide at the point $d.$
\end{theorem}

\begin{theorem}\label{main1 of example}Let $f(x)=(x^d,0,0,\cdots)+(\lambda x,0,0,\cdots)$, $p>2(d-1)^2+1.$ Then
$$T-\text{adic NP of }C_{f,u}(s,T)=\text{ord}_p(q)P_{\{1,d\},u}. $$
\end{theorem}
\begin{theorem}\label{main2 of example}
Let $f(x)=(x^d,0,0,\cdots)+(\lambda x,0,0,\cdots)$, $p>2(d-1)^2+1,$ and $m\geq1.$ Then
$$\pi_m-\text{adic NP of }C_{f,u}(s,\pi_m)=\text{ord}_p(q)P_{\{1,d\},u}.$$
\end{theorem}
\begin{theorem}\label{main3 of example}
Let $f(x)=(x^d,0,0,\cdots)+(\lambda x,0,0,\cdots)$, $p>2(d-1)^2+1,$ and $m\geq1.$ Then
$$\pi_m-\text{adic NP of }L_{f,u}(s,\pi_m)=\text{ord}_p(q)P_{\{1,d\},u}~~ on~~[0,p^{m-1}d].$$
\end{theorem}

The $q$-adic Newton polygon of $L_{f,0}(s,\pi_m)$ with $m=1$ was
considered by Roger Yang \cite{Ro} and Zhu \cite{Zh1}. In the first
paper, the author got all slopes of the Newton polygon for
$p\equiv-1(\mod d)$ and the first slope for all $p.$ In Zhu's paper,
she pointed that if $\lambda\neq0,$ the Newton polygon are all
generic, and the Newton polygon goes to Hodge as $p$ goes to
infinity, hence proved Wan's conjecture\cite{Wa2} for this certain
case. S. Sperber\cite{SS} also studied the Newton polygon of
$L_{f,0}(s,\pi_m)$ when $f$ is of degree 3 and $m=1$.

\section{Twisted T-adic Dwork's trace formula}
In this section we review the twisted $T$-adic analogy of Dwork's theory on exponential sums.

Let
$E(x)=\exp(\sum\limits_{i=0}^{\infty}\frac{x^{p^i}}{p^i})=\sum\limits_{n=0}^{\infty}\lambda_nx^n$
be the Artin-Hasse series. Define a new T-adic uniformizer $\pi$ of
the $T$-adic local ring $\mathbb{Q}_p[[T]]$ by the formula
$E(\pi)=1+T.$

Recall $C_u=C_u([0,d])=\{v\in\mathbb{N}|v\equiv u(\mod q-1)\}.$ Write
$$L_u=\{\sum_{v\in
C_u}b_v\pi^{\frac{v}{d(q-1)}}x^{\frac{v}{q-1}}:b_v\in\mathbb{Z}_q[[\pi^{\frac{1}{d(q-1)}}]]\}$$ and
$$B_u=\{\sum_{v\in
C_u}b_v\pi^{\frac{v}{d(q-1)}}x^{\frac{v}{q-1}}:b_v\in\mathbb{Z}_q[[\pi^{\frac{1}{d(q-1)}}]];{\rm
ord}_{\pi}b_v\rightarrow\infty,v\rightarrow\infty\}.$$

Define $\psi_p:L_u\rightarrow L_{p^{-1}u}$ by
$\psi_p(\sum\limits_{v\in C_u}b_vx^v)=\sum\limits_{v\in C_{p^{-1}u}}b_{pv}x^v.$ Then the map $\psi_p\circ E_f$
sends $L_u$ to $B_{p^{-1}u},$ where
$E_f(x)=E(\pi x^d)E(\pi\lambda x)=\sum\limits_{n=0}^{\infty}\gamma_nx^n.$

The Galois group ${\rm Gal}(\mathbb{Q}_q/\mathbb{Q}_p)$ can act on
$B_u$ by fixing $\pi^{\frac{1}{d(q-1)}}$ and $x.$ Let $\sigma\in{\rm
Gal}(\mathbb{Q}_q/\mathbb{Q}_p)$ be the Frobenius element such that
$\sigma(\zeta)=\zeta^p$ if $\zeta$ is a $(q-1)$-th root of unity.
The operator $\Psi=\sigma^{-1}\circ\psi_p\circ E_f(x)$ sends $B_u$
to $B_{p^{-1}u},$ hence $\Psi$ operators on
$B=\bigoplus\limits_{i=0}^{b-1}B_{p^iu}.$ We call it Dwork's
$T$-adic semi-linear operator because it is semi-linear over
$\mathbb{Z}_q[[\pi^{\frac{1}{d(q-1)}}]].$

We have $\Psi^a=\psi_p^a\circ
\prod_{i=0}^{a-1}E_f^{\sigma^i}(x^{p^i}).$ It follows that $\Psi^a$
operates on $B_u$ and is linear over
$\mathbb{Z}_p[[\pi^{\frac{1}{d(q-1)}}]].$ Moreover, one can show
that $\Psi$ is completely continuous in sense of Serre \cite{Se}, so
$\det(1-\Psi^as|B_u/\mathbb{Z}_q[[\pi^{\frac{1}{d(q-1)}}]])$ and
$\det(1-\Psi s|B/\mathbb{Z}_p[[\pi^{\frac{1}{d(q-1)}}]])$ are well
defined.

Now we state the twisted $T$-adic Dwork's trace formula \cite{Liu1}.
\begin{theorem}We have
$$C_{f,u}(s,T)=\det(1-\Psi^as|B_u/\mathbb{Z}_q[[\pi^{\frac{1}{d(q-1)}}]]).$$
\end{theorem}

\section{The twisted Dwork semi-linear operator}
In this section, we study the twisted Dwork semi-linear operator
$\Psi.$

Recall that $\Psi=\sigma^{-1}\circ\psi_p\circ E_f,$
$E_f(x)=E(\pi x^d)E(\pi\lambda x)=\sum\limits_{n=0}^{\infty}\gamma_nx^n,$
where $$\gamma_n=\sum_{di+j=n;i,j\geq0}\pi^{i+j}\lambda_i\lambda_j\hat{\lambda}^j.$$

We see $B=\bigoplus_{i=1}^bB_{p^iu}$ has a basis represented by
$\{x^{\frac{s_i}{q-1}+j}\}_{1\leq i\leq b,~j\in\mathbb{N}}$ over $\mathbb{Z}_q[[\pi^{\frac{1}{d(q-1)}}]].$
 We have
$$\Psi(x^{\frac{s_i}{q-1}+j})=\sigma^{-1}\circ\psi_p(\sum_{l=0}^{\infty}\gamma_lx^{\frac{s_i}{q-1}+j+l})
=\sum_{l=0}^{\infty}\gamma_{pl+u_{b-i}-j}^{\sigma^{-1}}x^{l+\frac{s_{i-1}}{q-1}}.$$

For $1\leq k,i\leq b$ and $l,j\in\mathbb{N},$ define
$$\gamma_{(\frac{s_k}{q-1}+l,\frac{s_i}{q-1}+j)}=\left\{
                        \begin{array}{ll}
                          \gamma_{pl+u_{b-i}-j}, & \hbox{} k=i-1;\\
                          0, & \hbox{} otherwise.
                        \end{array}
                      \right.
$$
Let
$\xi_1,\cdots,\xi_{a}$ be a normal basis of $\mathbb{Q}_q$ over
$\mathbb{Q}_p$ and write
$$\xi_v^{\sigma^{-1}}\gamma_{(\frac{s_k}{q-1}+l,\frac{s_i}{q-1}+j)}^{\sigma^{-1}}=
\sum_{w=1}^{a}\gamma_{(w,\frac{s_k}{q-1}+l)(v,\frac{s_i}{q-1}+j)}\xi_w.$$
It is easy to see $\gamma_{(w,\frac{s_k}{q-1}+l)(v,\frac{s_i}{q-1}+j)}=0$ for any $w$ and $v$ if $k\neq i-1.$
 Define the $i$-th submatrix by
$$\Gamma^{(i)}=(\gamma_{(w,\frac{s_{i-1}}{q-1}+l)(v,\frac{s_i}{q-1}+j)})_{1\leq w,v\leq a;l,j\in\mathbb{N}},$$
then the matrix of the operator $\Psi$ on B over
$\mathbb{Z}_p[[\pi^{\frac{1}{d(q-1)}}]]$ with respect to the basis
$\{\xi_vx^{\frac{s_i}{q-1}+j}\}_{1\leq i\leq b,1\leq v\leq a;j\in\mathbb{N}}$ is

$$\Gamma=\left(
                \begin{array}{ccccc}
                  0 & \Gamma^{(1)} & 0 & \cdots & 0 \\
                  0 & 0 & \Gamma^{(2)} & \cdots & 0 \\
                  \vdots & \vdots & \vdots & \vdots & \vdots \\
                  0 & 0 & 0 & \cdots & \Gamma^{(b-1)} \\
                  \Gamma^{(b)}& 0 & 0 & \cdots & 0 \\
                \end{array}
              \right).
$$

Hence by a
result of Li and Zhu \cite{LZ}, we have $$\det(1-\Psi
s|B/\mathbb{Z}_p[[\pi^{\frac{1}{d(q-1)}}]])=\det(1-\Gamma s)=\sum_{n=0}^{\infty}(-1)^{bn}C_{bn}s^{bn},$$
with $C_n=\sum\det(A)$ where
$A$ runs over all principle $n\times n$ submatrix  of $\Gamma.$

For every principle submatrix $A$ of $\Gamma,$ write
$A^{(i)}=A\cap\Gamma^{(i)}$ as the sub-matrix of $\Gamma^{(i)}.$ For
a principle $bn\times bn$ submatrix $A$ of $\Gamma,$ by linear
algebra, if one of $A^{(i)}$ is not $n\times n$ submatrix of
$\Gamma^{(i)}$ then at least one row or column of $A$ are 0 since
$A$ is principle.

Let $\mathcal{A}_n$ be the set of all $bn\times bn$ principle
submatrix A of $\Gamma$ with $A^{(i)}$ all $n\times n$ submatrix of
$\Gamma^{(i)}$ for each $1\leq i\leq b.$ Then we have
$$C_{bn}=\sum_{A\in\mathcal{A}_n}\det(A)=\sum_{A\in\mathcal{A}_n}(-1)^{n^2(b-1)}\prod_{i=1}^b\det(A^{(i)}).$$

Let $O(\pi^{\alpha})$ denotes any element of $\pi$-adic order $\geq
\alpha.$
\begin{lemma}We have
$$\gamma_n=\pi^{[\frac{n}{d}]+d\{\frac{n}{d}\}}\lambda_{[\frac{n}{d}]}\lambda_{d\{\frac{n}{d}\}}\hat{\lambda}^{d\{\frac{n}{d}\}}
+O(\pi^{[\frac{n}{d}]+d\{\frac{n}{d}\}+1}).$$
\end{lemma}
\proof
This follows from the fact that $i+j\geq[\frac{n}{d}]+d\{\frac{n}{d}\}$ if $di+j=n.$
\endproof
\begin{corollary}
For any $1\leq i\leq b$ and $1\leq w,u\leq a,$ we have
$$\gamma_{(w,\frac{s_{i-1}}{q-1}+l)(v,\frac{s_i}{q-1}+j)}=O(\pi^{[\frac{pl+u_{b-i}-j}{d}]+d\{\frac{pl+u_{b-i}-j}{d}\}}).$$
\end{corollary}

\proof By
$$\xi_v^{\sigma^{-1}}\gamma_{(\frac{s_{i-1}}{q-1}+l,\frac{s_i}{q-1}+j)}^{\sigma^{-1}}=\xi_v^{\sigma^{-1}}\gamma_{pl-j+u_{b-i}}^{\sigma^{-1}}
=\sum_{w=1}^{a}\gamma_{(w,\frac{s_{i-1}}{q-1}+l)(v,\frac{s_i}{q-1}+j)}\xi_w,$$ this follows from last
lemma.\endproof

\begin{theorem}\label{estimate for example}
Let
$R_i$ be a finite subset of $\{1,2,\cdots,a\}\times(\frac{s_i}{q-1}+\mathbb{N})$ of cardinality $an$ for each $1\leq i\leq b,$ $\tau$ a permutation of
$\bigcup_{i=1}^bR_i.$ If $p>2(d-1)^2+1,$ then
$$\sum_{i=1}^b\sum_{l\in R_i}([\frac{p\phi_i(l)+u_{b-i}-\phi_{i(\tau(l))}(\tau(l))}{d}]+d\{\frac{p\phi_i(l)+u_{b-i}-\phi_{i(\tau(l))}(\tau(l))}{d}\})$$
$$\geq abP_{\{1,d\},u}(n)+\sum_{i=1}^b\sharp\{l\in R_i|\phi_i(l)>n-1\},$$
where $i(l)=i$ if $l\in R_i$ and $\phi_i$ is the projection
$\{1,2,\cdots,a\}\times(\frac{s_i}{q-1}+\mathbb{N})\rightarrow \mathbb{N}$ such that $$\phi_i((\cdot,\frac{s_{i}}{q-1}+l))=l.$$
\end{theorem}
\proof
We have
$$\sum_{i=1}^b\sum_{l\in R_i}([\frac{p\phi_i(l)+u_{b-i}-\phi_{i(\tau(l))}(\tau(l))}{d}]+d\{\frac{p\phi_i(l)+u_{b-i}-\phi_{i(\tau(l))}(\tau(l))}{d}\})$$
$$=\sum_{i=1}^b\sum_{l\in R_i}\frac{(p-1)\phi_i(l)+u_{b-i}}{d}+(d-1)\sum_{i=1}^b\sum_{l\in R_i}(\{\frac{p\phi_i(l)+u_{b-i}}{d}\}-\{\frac{\phi_i(l)}{d}\})$$
$$+(d-1)\sum_{i=1}^b\sum_{l\in R_i}1_{\{\frac{p\phi_i(l)+u_{b-i}}{d}\}<\{\frac{\phi_{i(\tau(l))}(\tau(l))}{d}\}}.$$
Observe that
$\sharp\{l\in R_i|\phi_i(l)>n-1\}=\sharp\{l\not\in R_i|\phi_i(l)\leq n-1\},$ we have
$$\sum_{l\in R_i}\frac{(p-1)\phi_i(l)+u_{b-i}}{d}$$
$$=(\sum_{l\in\phi_i^{-1}(A_n)}+\sum_{l\in R_i,\phi_i(l)>n-1}+\sum_{\phi_i(l)\leq n-1,l\not\in R_i})\frac{(p-1)\phi_i(l)+u_{b-i}}{d}$$
$$\geq a\sum_{l=0}^{n-1}\frac{(p-1)l+u_{b-i}}{d}+\frac{p-1}{d}\sharp\{l\in R_i|\phi_i(l)>n-1\}.$$
Similarly, dividing the sum into three disjoint parts as above, we have
$$\sum_{l\in R_i}\{\frac{p\phi_i(l)+u_{b-i}}{d}\}-\{\frac{\phi_i(l)}{d}\}$$
$$\geq a\sum_{l=0}^{n-1}(\{\frac{pl+u_{b-i}}{d}\}-\{\frac{l}{d}\})-\frac{2(d-1)^2}{d}\sharp\{l\in R_i|\phi_i(l)>n-1\}.$$
The theorem follows from $p>2(d-1)^2+1.$
\endproof

\begin{theorem}\label{order of C_bn}
If $p>2(d-1)^2+1,$ then we have $$\text{ord}_{\pi}(C_{abn})\geq
abP_{\{1,d\},u}(n).$$ In particular, $$C_{abn}=\pm{\rm
Norm}(\prod_{i=1}^b\det((\gamma_{pl-j+u_{b-i}})_{l,j\in
A_n}))+O(\pi^{abP_{\{1,d\},u}(n)+\frac{1}{d(q-1)}}),$$
where Norm is the norm map from
$\mathbb{Q}_q(\pi^{\frac{1}{d(q-1)}})$ to
$\mathbb{Q}_p(\pi^{\frac{1}{d(q-1)}}).$
\end{theorem}

\proof Let $A\in \mathcal{A}_{an},$ $R_i$ the set of rows of
$A^{(i)}$ as the submatrix of $\Gamma^{(i)},$ $\tau$ a permutation
of $\bigcup_iR_i.$ By the above corollary and Theorem
\ref{estimate for example}, we have
$$\text{ord}_{\pi}(C_{abn})\geq\sum_{i=1}^b\sum_{l\in R_i}\lceil\frac{p\phi_i(l)+u_{b-i}-\phi_{i(\tau(l))}(\tau(l))}{d}\rceil
\geq abP_{\{1,d\},u}(n).$$ Moreover the strict inequality holds if
there exist $1\leq i\leq b$ such that $R_i\neq\phi_i^{-1}(A_n),$ hence
$$C_{abn}=\prod_{i=1}^b\det(\gamma_{(w,\frac{s_{i-1}}{q-1}+l)(v,\frac{s_i}{q-1}+j)})_{1\leq w,v\leq a;l,j\in A_n}+O(\pi^{abP_{\{1,d\},u}(n)+\frac{1}{d(q-1)}}).$$
Therefore the theorem follows from the following.
\endproof

\begin{lemma}\label{dimup}
For finite subset $A\subset\mathbb{N},$ we have
$$\det((\gamma_{(w,\frac{s_{i-1}}{q-1}+l)(v,\frac{s_i}{q-1}+j)})_{1\leq w,v\leq a;l,j\in A})=\pm{\rm Norm}(\det(\gamma_{(\frac{s_{i-1}}{q-1}+l,\frac{s_i}{q-1}+j)})_{l,j\in A})),$$
where Norm is the norm map from
$\mathbb{Q}_q(\pi^{\frac{1}{d(q-1)}})$ to
$\mathbb{Q}_p(\pi^{\frac{1}{d(q-1)}}).$
\end{lemma}
\proof Let $V=\oplus_{j\in
A}\mathbb{Q}_q(\pi^{\frac{1}{d(q-1)}})e_j$ be a vector space over
$\mathbb{Q}_q(\pi^{\frac{1}{d(q-1)}})$, and let $F$ be the linear
operator on $V$ whose matrix with respect to the basis $\{e_j\}$ is
$(\gamma_{(\frac{s_{i-1}}{q-1}+l,\frac{s_i}{q-1}+j)})_{l,j\in A}$, and let $\sigma$ act on V
coordinate-wise. It is easy to see that
$(\gamma_{(w,\frac{s_{i-1}}{q-1}+l)(v,\frac{s_i}{q-1}+j)})_{1\leq w,v\leq a;l,j\in A}$ is the matrix of
$\sigma^{-1}\circ F$ over $\mathbb{Q}_p(\pi^{\frac{1}{d(q-1)}})$
with respect to the basis $\{\xi_ue_j\}_{1\leq u\leq a;l\in A}$.
Therefore
$$\det((\gamma_{(w,\frac{s_{i-1}}{q-1}+l)(v,\frac{s_i}{q-1}+j)})_{1\leq w,v\leq a;l,j\in A})=\det(\sigma^{-1}\circ F)$$$$=\pm\det(F\mid
\mathbb{Q}_p(\pi^{\frac{1}{d(q-1)}}))=\pm{\rm Norm}(\det(F)).$$ The
lemma is proved.\endproof

\section{Hasse polynomial}
In this section, $1\leq n\leq d-1.$ We shall study
$\det((\gamma_{pl-j+u_{b-i}})_{0\leq l,j\leq n-1})).$
\begin{definition}
For any $1\leq n\leq d-1$ and $1\leq i\leq b,$ we define
$$S_{n,i}=\{\tau\in S_n|\tau(l)\leq d\{\frac{pl+u_{b-i}}{d}\}\},$$
$$H_{n,u}^{(i)}(y)=\sum_{\tau\in S_{n,i}}\text{sgn}(\tau)\prod_{l=0}^{n-1}\lambda_{[\frac{pl+u_{b-i}-\tau(l)}{d}]}\lambda_{d\{\frac{pl+u_{b-i}-\tau(l)}{d}\}}y^{d\{\frac{pl+u_{b-i}-\tau(l)}{d}\}}.$$
The $u$-twisted Hasse polynomial $H_{n,u}(y)$ of $\{1,d\}$ at $n$ is defined by
$$H_{n,u}(y)=\prod_{i=1}^bH_{n,u}^{(i)}(y)(\mod p).$$
The $u$-twisted Hasse polynomial $H_u(y)$ of $\{1,d\}$ is defined by
$H_u=\prod\limits_{n=1}^{d-1}H_{n,u}.$
\end{definition}

\begin{lemma}\label{momomial}
We have $$H_{n,u}^{(i)}(y)=y^{\sum\limits_{l=0}^{n-1}(d\{\frac{pl+u_{b-i}}{d}\}-l)}\sum_{\tau\in S_{n,i}}\text{sgn}(\tau)\prod_{l=0}^{n-1}
\lambda_{[\frac{pl+u_{b-i}}{d}]}\lambda_{d\{\frac{pl+u_{b-i}}{d}\}-\tau(l)}.$$
\end{lemma}
\proof The lemma follows from
$$\sum_{l=0}^{n-1}d\{\frac{pl+u_{b-i}-\tau(l)}{d}\}=\sum_{l=0}^{n-1}(d\{\frac{pl+u_{b-i}}{d}\}-l),$$
for any $\tau\in S_{n,i}.$
\endproof

\begin{theorem}
For $1\leq n\leq d-1$ and $p>2(d-1)^2+1,$ we have
$$C_{abn}=\pm\text{Norm}(H_{n,u}(\hat{\lambda}))\pi^{abP_{\{1,d\},u}(n)}+O(\pi^{abP_{\{1,d\},u}(n)+\frac{1}{d(q-1)}}).$$
\end{theorem}

\proof  We have
 $$\det(\gamma_{pl-j+u_{b-i}})_{0\leq l,j\leq
n-1}=\sum_{\tau\in S_n}{\rm
sgn}(\tau)\prod_{l=0}^{n-1}\gamma_{pl-\tau(l)+u_{b-i}}.$$
For any $\tau\in S_n,$ we have
$$\sum_{l=0}^{n-1}\text{ord}_{\pi}\gamma_{pl+u_{b-i}-\tau(l)}
\geq\sum_{l=0}^{n-1}\frac{pl+u_{b-i}-\tau(l)}{d}+(d-1)\{\frac{pl+u_{b-i}-\tau(l)}{d}\}$$
$$\geq\sum_{l=0}^{n-1}\frac{(p-1)l+u_{b-i}}{d}+(d-1)\sum_{l=0}^{n-1}(\{\frac{pl+u_{b-i}}{d}\}-\{\frac{l}{d}\}+1_{\tau(l)>d\{\frac{pl+u_{b-i}}{d}\}})$$
$$\geq\sum_{l=0}^{n-1}\frac{(p-1)l+u_{b-i}}{d}+(d-1)\sum_{l=0}^{n-1}(\{\frac{pl+u_{b-i}}{d}\}-\{\frac{l}{d}\}),$$
with the equalities holding if and only if $\tau\in S_{n,i}.$ The
lemma follows.
\endproof

\begin{theorem}
The Hasse polynomial $H_u$ is nonzero.
\end{theorem}

\proof By Lemma \ref{momomial}, we see $H_{n,u}^{(i)}$ is a
monomial, so is $H_u$. It suffices to show for any $1\leq i\leq b$
and $1\leq n\leq d,$ the coefficients of $H_{n,u}^{(i)}$
$$f_{n,p}^{(i)}=\sum_{\tau\in S_{n,i}}\text{sgn}(\tau)\prod_{l=0}^{n-1}\lambda_{[\frac{pl+u_{b-i}}{d}]}\lambda_{d\{\frac{pl+u_{b-i}}{d}\}-\tau(l)}\in\mathbb{Z}_p^{\times}.$$
Write
$u_{n,i}=\prod_{l=0}^{n-1}[\frac{pl+u_{b-i}}{d}]!(d\{\frac{pl+u_{b-i}}{d}\})!\in\mathbb{Z}_p^{\times}$
and $\alpha_l=d\{\frac{pl+u_{b-i}}{d}\},$ it suffices to show
$u_{n,i}f_{n,p}^{(i)}\in\mathbb{Z}_p^{\times}.$
$$u_{n,i}f_{n,p}^{(i)}
=\sum_{\tau\in S_{n,i}}\text{sgn}(\tau)\prod_{l=0}^{n-1}\alpha_l(\alpha_l-1)\cdots(\alpha_l-\tau(l)+1)$$
$$=\sum_{\tau\in S_n}\text{sgn}(\tau)\prod_{l=0}^{n-1}\alpha_l(\alpha_l-1)\cdots(\alpha_l-\tau(l)+1)\\
=\det M,$$
where $$M=\left(
           \begin{array}{ccccc}
             1 & \alpha_0 & \alpha_0(\alpha_0-1) & \cdots & \alpha_0\cdots(\alpha_0-n+2) \\
             1 & \alpha_1 & \alpha_1(\alpha_1-1) & \cdots & \alpha_1\cdots(\alpha_1-n+2)\\
             \vdots& \vdots & \vdots & \vdots & \vdots \\
             1 & \alpha_{n-1} &\alpha_{n-1} (\alpha_{n-1}-1) & \cdots & \alpha_{n-1}\cdots(\alpha_{n-1}-n+2) \\
           \end{array}
         \right).
$$
It is easy to see that $M$ is equivalent to a Vandermonde matrix. So
$$u_{n,i}f_{n,p}^{(i)}=\det M=\prod\limits_{0\leq l<j\leq
n-1}(\alpha_l-\alpha_j)\in\mathbb{Z}_p^{\times},$$ the theorem is
proved.
\endproof

\begin{theorem}\label{independency of lambda}
For $1\leq n\leq d-1$ and $p>2(d-1)^2+1,$ we have
$$\text{ord}_{\pi}(C_{abn})=abP_{\{1,d\},u}(n).$$
\end{theorem}
\proof
The theorem follows from the last two theorems and Lemma \ref{momomial}.
\endproof

\section{Proof of the main theorems}
In this section we prove the main theorems of this paper.

Firstly we prove Theorem \ref{main of example}, which says that $$P_{\{1,d\},u}\geq (p-1)H_{[0,d],u}^{\infty}$$
with equality holding at the point $d.$

{\it Proof of Theorem \ref{main of example}.} It need only to show this for $n\leq d.$
 The inequality follows from the following
$$\sum_{l=0}^{n-1}\{\frac{pl+u_i}{d}\}\geq\sum_{l=0}^{n-1}\{\frac{l}{d}\}.$$
Moreover, the equality above holds if $n=d,$ the theorem follows.\qed

\begin{lemma}\label{zqtozp}
The $T$-adic Newton polygon of
$\det(1-\Psi^as^a|B/\mathbb{Z}_q[[\pi^{\frac{1}{d(q-1)}}]])$
coincides with that of $\det(1-\Psi
s|B/\mathbb{Z}_p[[\pi^{\frac{1}{d(q-1)}}]])$
\end{lemma}

\proof The lemma follows from the following:
\begin{align*}
\prod_{\zeta^a=1}\det(1-\Psi\zeta
s|B/\mathbb{Z}_p[[\pi^{\frac{1}{d(q-1)}}]])&=\det(1-\Psi^as^a|B/\mathbb{Z}_p[[\pi^{\frac{1}{d(q-1)}}]])\\
&={\rm
Norm}(\det(1-\Psi^as^a|B/\mathbb{Z}_q[[\pi^{\frac{1}{d(q-1)}}]])),
\end{align*}
where Norm is the norm map from
$\mathbb{Q}_q[[\pi^{\frac{1}{d(q-1)}}]]$ to
$\mathbb{Q}_p[[\pi^{\frac{1}{d(q-1)}}]].$
\endproof

\begin{lemma}\label{BtoB_i}
The T-adic Newton polygon of $C_{f,u}(s,T)^b$ coincides with that of
$\det(1-\Psi^as|B/\mathbb{Z}_q[[\pi^{\frac{1}{d(q-1)}}]]).$
\end{lemma}

\proof Let $\text{Gal}(\mathbb{Q}_q/\mathbb{Q}_p)$ act on $\mathbb{Z}_q[[T]][[s]]$ by fixing $s$ and $T,$
then we have $$C_{pu,f}(s,T)=C_{f,u}(s,T)^{\sigma}.$$
Therefore the lemma follows from the following
$$\prod_{i=0}^{b-1}C_{f,u}(s,T)^{\sigma^i}=\prod_{i=0}^{b-1}C_{p^iu,f}(s,T)$$
$$=\prod_{i=0}^{b-1}\det(1-\Psi^as|B_{p^iu}/\mathbb{Z}_q[[\pi^{\frac{1}{d(q-1)}}]])
=\det(1-\Psi^as|B/\mathbb{Z}_q[[\pi^{\frac{1}{d(q-1)}}]]).$$
\endproof

\begin{theorem}\label{description of NP}
The $T$-adic Newton polygon of $C_{f,u}(s,T)$ is the lower convex closure of the points
$$(n,\frac{1}{b}\text{ord}_{T}(C_{abn})),~n=0,1,\cdots.$$
\end{theorem}
\proof By Lemma \ref{zqtozp}, we see the T-adic Newton polygon of
the power series
$\det(1-\Psi^as^a|B/\mathbb{Z}_q[[\pi^{\frac{1}{d(q-1)}}]])$ is the
lower convex closure of the points $$(n,{\rm
ord}_{T}C_{n}),~~~n=0,1,\cdots.$$ It is clear that
$(n,\text{ord}_{T}(c_n))$ is not a vertex of that polygon if $a\nmid
n$. So that Newton polygon is the lower convex closure of the points
$$(an,\text{ord}_{T}(c_{an})),\ n=0,1,\cdots.$$
Hence the T-adic Newton polygon of
$\det(1-\Psi^as|B/\mathbb{Z}_q[[\pi^{\frac{1}{d(q-1)}}]])$ is the
convex closure of the points $$(n,{\rm
ord}_{T}C_{an}),~~~n=0,1,\cdots.$$ By Lemma \ref{BtoB_i}, the T-adic
Newton polygon of $C_{f,u}(s,T)^b$ is the lower convex closure of
the points $$(n,{\rm ord}_{T}C_{an}),~~~n=0,1,\cdots,$$ hence the
closure of the points $$(bn,{\rm ord}_{T}C_{abn}),~~~n=0,1,\cdots.$$
It follows that the $T$-adic Newton polygon of $C_{f,u}(s,T)$ is the
convex closure of the points
$$(n,\frac{1}{b}\text{ord}_{T}(C_{abn})),~n=0,1,\cdots.$$ The theorem is proved.
\endproof
\begin{lemma}\label{geq}Let $f(x)=(x^d,0,0,\cdots)+(\lambda x,0,0,\cdots)$, $p>2(d-1)^2+1.$
Then $$T-\text{adic NP of }C_{f,u}(s,T)\geq\text{ord}_p(q)P_{\{1,d\},u}. $$
\end{lemma}
\proof By the last theorem, we see that
the $T$-adic Newton polygon of $C_{f,u}(s,T)$ is the convex closure of the points
$$(n,\frac{1}{b}\text{ord}_{T}(C_{abn})),~n=0,1,\cdots.$$ The theorem now follows from Theorem \ref{order of C_bn}.
\endproof

\begin{theorem}\label{spe} Let $A(s,T)$ be a
$T$-adic entrie series in $s$ with unitary constant term. If
$0\neq|t|_p<1$, then
$$t-adic\text{ NP of
}A(s,t)\geq T-adic\text{ NP of }A(s,T),$$ where NP is the short for
Newton polygon. Moreover, the equality holds for one $t$ if and only
if it holds for all $t$.\end{theorem}
\proof
The reader may refer \cite{LLN} and we omit the proof here.
\endproof

\begin{theorem}Let $f(x)=(x^d,0,0,\cdots)+(\lambda x,0,0,\cdots)$, $p>2(d-1)^2+1.$
If the equality $$\pi_m-\text{adic NP of }C_{f,u}(s,\pi_m)={\rm
ord}_p(q)P_{\{1,d\},u}$$ holds for one $m\geq1$, then it holds for all
$m\geq1$, and we have $$T-\text{adic NP of }C_{f,u}(s,T)={\rm
ord}_p(q)P_{\{1,d\},u}.$$
\end{theorem}
\proof
It follows from Lemma \ref{geq} and the last theorem.
\endproof

\begin{theorem}Let $f(x)=(x^d,0,0,\cdots)+(\lambda x,0,0,\cdots)$, $p>2(d-1)^2+1.$
Then  $$\pi_m-\text{adic NP of }C_{f,u}(s,\pi_m)={\rm
ord}_p(q)P_{\{1,d\},u}$$ if and only if
$$\pi_m-\text{adic NP of }L_{f,u}(s,\pi_m)={\rm
ord}_p(q)P_{\{1,d\},u}\text{ on }[0,p^{m-1}d].$$
\end{theorem}
\proof Assume that $ L_{f,u}(s,
\pi_m)=\prod\limits_{i=1}^{p^{m-1}d}(1-\beta_is) $. Then
$$C_{f,u}(s,\pi_m )= \prod\limits_{j=0}^{\infty} L_{f,u}(q^js, \pi_m)=
\prod\limits_{j=0}^{\infty}\prod\limits_{i=1}^{p^{m-1}d}(1-\beta_iq^js).$$
Therefore the slopes of the  $q$-adic Newton polygon of
$C_{f,u}(s,\pi_m)$ are the numbers
$$j+{\rm ord}_q(\beta_i),\ 1\leq
i\leq p^{m-1}{\rm Vol}(\triangle), j=0,1,\cdots.$$
Since $$w(n+m^{m-1}d)=p^m-p^{m-1}+w(n),$$
then the slopes of $P_{\{1,d\},u}$ are the numbers
 $$j(p^m-p^{m-1})+w(i),\ 1\leq i\leq
p^{m-1}d, j=0,1,\cdots.$$ It follows
that
$$\pi_m-\text{adic NP of }C_{f,u}(s,\pi_m)={\rm
ord}_p(q)P_{\{1,d\},u}$$ if and only if
$$\pi_m-\text{adic NP of }L_{f,u}(s,\pi_m)={\rm
ord}_p(q)P_{\{1,d\},u}\text{ on }[0,p^{m-1}d].$$
\endproof

We now prove Theorems \ref{main1 of example}, \ref{main2 of example} and
\ref{main3 of example}. By
the above theorems, it suffices to prove the following.
\begin{theorem}\label{final th}
Let $f(x)=(x^d,0,0,\cdots)+(\lambda x,0,0,\cdots)$, $p>2(d-1)^2+1.$
Then  $$\pi_1-\text{adic NP of }L_{f,u}(s,\pi_1)={\rm
ord}_p(q)P_{\{1,d\},u}\text{ on }[0,d]$$ if and
only if $H((a_i)_{0\leq i\leq d})\neq0$.
\end{theorem}\proof By a result of Liu \cite{Liu2}, the
$q$-adic Newton polygon  of $L_{f,u}(s,\pi_1)$ coincides with
$H_{[0,d],u}^{\infty}$ at the point $d$. By
Theorem \ref{main of example}, $P_{\{1,d\},u}$ coincide with
$(p-1)H_{[0,d],u}^{\infty}$ at the point $d$,
It follows that the $\pi_1$-adic Newton polygon of $L_{f,u}(s,\pi_1)$
coincides with ${\rm ord}_p(q)p_{\triangle}$ at the point $d$. Therefore it suffices to show that
$$\pi_1-\text{adic NP of }L_{f,u}(s,\pi_1)={\rm
ord}_p(q)P_{\{1,d\},u}\text{ on }[0,d-1]$$ if and
only if $H((a_u)_{u\in\triangle})\neq0$.

From the identity $$C_{f,u}(s,\pi_1 )= \prod\limits_{j=0}^{\infty}
L_{f,u}(q^js, \pi_1),$$ and the fact the $q$-adic orders of the
reciprocal roots of $L_{f,u}(s,\pi_1)$ are no greater than $1$, we
infer that
$$\pi_1-\text{adic NP of }L_{f,u}(s,\pi_1) =\pi_1-\text{adic NP of
}C_{f,u}(s,\pi_1)\text{ on }[0,d-1].$$Therefore it
suffices to show that $$\pi_1-\text{adic NP of }C_{f,u}(s,\pi_1)={\rm
ord}_p(q)P_{\{1,d\},u}\text{ on }[0,d-1]$$ if and
only if $H((a_i)_{0\leq i\leq d})\neq0$. The theorem now follows
from the $T$-adic Dwork trace formula and Theorems
\ref{independency of lambda}.\endproof

%=========================================================

\end{document}